\documentclass[12pt, leqno]{amsart}
\usepackage[mathscr]{eucal}
\usepackage{amssymb,amsmath}
\usepackage{amsfonts}
\usepackage{a4}

\newtheorem{theorem}{Theorem}
\newtheorem{proposition}{Proposition}
\newtheorem{lemma}{Lemma}
\newtheorem{corollary}{Corollary}

\begin{document}
\title{Deformations of monoidal functors}

\author[Tomasz Maszczyk]{Tomasz Maszczyk\dag}
\address{Institute of Mathematics\\
Polish Academy of Sciences\\
Sniadeckich 8\newline 00--956 Warszawa, Poland\\
\newline Institute of Mathematics\\
University of Warsaw\\ Banacha 2\newline 02--097 Warszawa, Poland}
\email{maszczyk@mimuw.edu.pl}

\thanks{\dag The author was partially supported by the grant 1261/7.PR UE/2009/7.}

\begin{abstract} We point out that for Yetter's deformational Hochschild  complex of a monoidal functor between abelian monoidal categories the Gerstenhaber-Voronov type operations can be defined  making it a strong homotopy Gerstenhaber algebra. This encodes deformation theory of monoidal functors in an analogical way as deformation theory of associative algebras is described by the  strong homotopy Gerstenhaber algebra structure on the corresponding Hochschild cochains. We describe a quasiclassical limit of deformations of a symmetric monoidal functor in terms of Poisson type structure. 
\end{abstract}

\maketitle

\subsection{Introduction}
The Gerstenhaber-Voronov operations \cite{ger-vor} lift to an action of an operad equivalent to the chain little discs operad on Tamarkin's resolution of the Gerstenhaber operad in the model structure on operads \cite{Tam}. This provided the first proof of Deligne's conjecture: the Gerstenhaber algebra structure on the Hochschild cohomology lifts to an algebra over an operad equivalent to the chain little discs operad structure on Hochschild cochains of an associative algebra. According to Yetter \cite{Yett,Yett'}  the most natural generalization of the construction of the Hochschild cochain complex for an associative algebra is the Hochschild cochain complex of a monoidal functor between abelian monoidal categories. As for associative algebras, 
first order deformations of monoidal functors are classified by the second Hochschild cohomology, while obstructions to higher order deformations lie in the third cohomology.
As Yetter showed, also the Gerstenhaber algebra structure on the Hochschild cohomology generalizes to monoidal functors as a result of generalization of the pre-Lie structure defined on Hochschild cochains by Gerstenhaber \cite{Ger}. 

The main point of the present paper consists in pointing out that the Gerstenha\break
ber-Voronov operations can be generalized in an analogical way from associative algebras to monoidal functors. According to Tamarkin's construction this implies a version of Deligne's conjecture for the Hochschild cochains of a monoidal functor.

Following Yetter's philosophy we point out that symmetric monoidal functors should play among monoidal functors an analogical role as commutative algebras play among associative ones. Indeed, we show that the  quasi-classical limit of a formal monoidal deformation of a symmetric monoidal functor still  makes sense and defines a kind of Poisson structure on a given symmetric monoidal functor. On the other hand, we show that natural transformations evaluated at the monoidal unit lead to the deformation theory of algebras (aka monoids in abelian monoidal categories). This suggests that monoidal functors could be considered as a wider foundation for noncommutative (non-symmetric monoidal) geometry than associative algebras.

\subsection{Monoidal functors} We will think about abelian monoidal categories as if there were categories of quasi-coherent sheaves on quasi-compact quasi-separated schemes. Every morphism $f:X\rightarrow
Y$ of such schemes induces a pair of adjoint functors $f^{*}\dashv f_{*}:{\rm Qcoh}_{X}\rightarrow
{\rm Qcoh}_{Y}$ (inverse and direct image) where the right adjoint (direct image) is monoidal. Motivated by this special case we borrow traditional notation from algebraic geometry and try to imagine a morphism of imaginary spaces, whose quasi-coherent direct image is our  monoidal functor. In this way deformations of monoidal functors get some geometric flavor, providing guiding principles in an abstract case. To make the theory easier we abuse language of monoidal categories on the ground of the Mac Lane coherence. We will also use possibly large groups or rings (with possibly proper class of elements) of natural transformations, since we do not restrict ourselves to small categories. This should not however cause troubles because all algebraic computations are performed object-wise. To avoid massive formulas or diagrams sometimes we use the element-wise convention as we were in a concrete category, assuming that the abstract meaning of this concrete shortcut is obvious.  

The monoidal structure of $f_{*}:{\rm Qcoh}_{X}\rightarrow
{\rm Qcoh}_{Y} $ consists of a morphism
\begin{align}\eta: {\mathcal
O}_{Y}\rightarrow f_{*}
{\mathcal O}_{X}\label{mon1}\end{align}
and a natural transformation $\mu$ of bifunctors
\begin{align}\mu^{{\mathcal F}_{1},{\mathcal F}_{2}}: f_{*}{\mathcal
F}_{1}\otimes f_{*}{\mathcal F}_{2}\rightarrow f_{*}({\mathcal
F}_{1}\otimes {\mathcal F}_{2})\label{mon2}\end{align}
making  the appropriate diagrams commute, which is equivalent to the following identities
\begin{align}\mu^{{\mathcal F},{\mathcal O}_{X}}(Id^{f_{*}{\mathcal
F}}\otimes \eta)=Id^{f_{*}{\mathcal F}},\ \ \mu^{{\mathcal O}_{X},{\mathcal F}}(\eta\otimes Id^{f_{*}{\mathcal
F}})=Id^{f_{*}{\mathcal F}},\label{un}\end{align}
\begin{align}\mu^{{\mathcal F}_{1},{\mathcal F}_{2}\otimes{\mathcal F}_{3}}
(Id^{f_{*}{\mathcal F}_{1}}\otimes \mu^{{\mathcal F}_{2},{\mathcal F}_{3}})
=\mu^{{\mathcal F}_{1}\otimes{\mathcal F}_{2},{\mathcal F}_{3}}
(\mu^{{\mathcal F}_{1},{\mathcal F}_{2}}
\otimes Id^{f_{*}{\mathcal F}_{3}})\label{ass}
\end{align}

\subsection{Deformations}
\paragraph{\textbf{Definition.}} We define a \emph{formal deformation of the monoidal structure} $(\eta, \mu)$ as a
sequence  $(\mu_{0}, \mu_{1}, \ldots )$ of natural transformations
\begin{align}\mu_{i}^{{\mathcal F}_{1},{\mathcal F}_{2}}: f_{*}{\mathcal
F}_{1}\otimes f_{*}{\mathcal F}_{2}\rightarrow f_{*}({\mathcal
F}_{1}\otimes {\mathcal F}_{2})\label{defmon2}\end{align}
such that $\mu_{0}=\mu$ and for $k>0$ the following identities hold
\begin{align}\sum_{i+j=k}\mu_{i}^{{\mathcal F}_{1},{\mathcal F}_{2}\otimes{\mathcal F}_{3}}
(Id^{f_{*}{\mathcal F}_{1}}\otimes \mu_{j}^{{\mathcal F}_{2},{\mathcal F}_{3}})
=\sum_{i+j=k}\mu_{i}^{{\mathcal F}_{1}\otimes{\mathcal F}_{2},{\mathcal F}_{3}}
(\mu_{j}^{{\mathcal F}_{1},{\mathcal F}_{2}}
\otimes Id^{f_{*}{\mathcal F}_{3}}).\label{defass}
\end{align}
Taking sequences of length $n+1$ $(\mu_{0}, \ldots,\mu_{n})$ such that $\mu_{0}=\mu$ and (\ref{defass}) holds for $k=1,\ldots,n$ we obtain the definition of an \emph{n-th infinitesimal deformation of the monoidal structure} $(\eta, \mu)$.

\vspace{3mm}
\paragraph{\textbf{Definition.}} We define the group of sequences $\varphi=(\varphi_{0}, \varphi_{1},  \ldots)$ of natural transformations
\begin{align}
\varphi_{i}^{{\mathcal F}}: f_{*}{\mathcal F}\rightarrow f_{*}{\mathcal F}
\end{align}
such that $\varphi_{0}^{{\mathcal F}}=Id^{f_{*}{\mathcal F}_{3}}$, with the neutral element $(Id, 0, 0, \ldots)$ and the composition
\begin{align}
(\varphi\widetilde{\varphi})_{k}^{{\mathcal F}}=\sum_{i+j=k}\varphi_{i}^{{\mathcal F}}\widetilde{\varphi}_{j}^{{\mathcal F}}.
\end{align}
We say that formal deformations $(\mu_{0}, \mu_{1}, \ldots )$ and $(\widetilde{\mu}_{0}, \widetilde{\mu}_{1}, \ldots )$ are \emph{equivalent} if there exists a sequence $(\varphi_{0}, \varphi_{1}, \ldots)$ such that for $k>0$
\begin{align}
\sum_{i+j=k}\varphi_{i}^{{\mathcal F}_{1}\otimes{\mathcal F}_{2}}\mu_{j}^{{\mathcal F}_{1},{\mathcal F}_{2}}=\sum_{i+j_{1}+j_{2}=k}\widetilde{\mu}_{i}^{{\mathcal F}_{1},{\mathcal F}_{2}}(\varphi_{j_{1}}^{{\mathcal F}_{1}}\otimes\varphi_{j_{2}}^{{\mathcal F}_{2}}).\label{equivass}
\end{align}
Taking sequences of length two $(\varphi_{0}, \varphi_{1})$ and 1-st infinitesimal deformations $(\mu_{0}, \mu_{1})$ and $(\widetilde{\mu}_{0}, \widetilde{\mu}_{1})$ such that $\varphi_{0}=Id$ and (\ref{equivass}) holds for $k=1$ we obtain the definition of  \emph{equivalence of 1-st infinitesimal deformations}.

\subsection{Hochschild complex}
\paragraph{\textbf{Definition.}} We define abelian groups of
$k$-cochains as follows.

For $k=0$ it consists of morphisms
\begin{align}
c^{0}: {\mathcal
O}_{Y}\rightarrow f_{*}
{\mathcal O}_{X}.\label{cochain1}
\end{align}

For $k>0$ it consists of natural transformations $c^{k}$ of
multifunctors
\begin{align}c^{{\mathcal F}_{1},\ldots,{\mathcal F}_{k}}: f_{*}{\mathcal
F}_{1}\otimes\cdots\otimes f_{*}{\mathcal F}_{k}\rightarrow
f_{*}({\mathcal F}_{1}\otimes\cdots\otimes {\mathcal
F}_{k}).\label{cochain2}\end{align}
In the element-wise convention
\begin{align}1 & \mapsto c^{0}(1)\end{align}
and for $k>0$
\begin{align}n_{1}\otimes\cdots\otimes n_{k} & \mapsto
c^{{\mathcal F}_{1},\ldots,{\mathcal F}_{k}}(n_{1},\ldots,n_{k}).\end{align}
We equip them with the Hochschild type differential \cite{hoch} as follows.
\begin{align}
&({\rm d}c^{0})^{{\mathcal F}_{1}}(n_{1}):=\mu^{{\mathcal F}_{1},{\mathcal O}_{X}}(n_{1},c^{0}(1))-\mu^{{\mathcal O}_{X}, {\mathcal F}_{1}}(c^{0}(1),n_{1})\label{diff1}
\end{align}
and for $k>0$
\begin{align}
&({\rm d}c^{k})^{{\mathcal F}_{1},\ldots,{\mathcal F}_{k+1}}(n_{1},\ldots,n_{k+1})\label{diff2}\\
&:= (-1)^{k-1}\mu^{{\mathcal F}_{1},{\mathcal F}_{2}\otimes\cdots\otimes{\mathcal F}_{k}}(n_{1}, c^{{\mathcal F}_{2},\ldots,{\mathcal F}_{k+1}}(n_{2},\ldots,n_{k+1}))\nonumber\\
&+\sum_{i=1}^{k}(-1)^{i+k-1}c^{{\mathcal F}_{1},\ldots,{\mathcal F}_{i}\otimes{\mathcal F}_{i+1},\ldots,{\mathcal F}_{k+1}}(n_{1},\ldots,\mu^{{\mathcal F}_{i},{\mathcal F}_{i+1}}(n_{i},n_{i+1}),\ldots, n_{k+1})\nonumber\\
&+\mu^{{\mathcal F}_{1}\otimes\cdots\otimes{\mathcal F}_{k},{\mathcal F}_{k+1}}(c^{{\mathcal F}_{1},\ldots,{\mathcal F}_{k}}(n_{1},\ldots,n_{k}),n_{k+1}).\nonumber
\end{align}

By a routine checking we see that $d^{2}=0$. We call the resulting cohomology \emph{Hochschild cohomology of a monoidal functor}.

We define the Gerstenhaber-Voronov type operations \cite{ger}\cite{ger-vor} on the cochains
\begin{align}
&(c^{i} c^{j})^{{\mathcal F}_{1},\ldots,{\mathcal F}_{i+j}}(n_{1},\ldots,n_{i+j})\label{prod}\\
&:= (-1)^{ij}\mu^{{\mathcal F}_{1}\otimes\cdots\otimes{\mathcal F}_{i}, {\mathcal F}_{i+1}\otimes\cdots\otimes{\mathcal F}_{i+j}}(c^{{\mathcal F}_{1},\ldots,{\mathcal F}_{i}}(n_{1},\ldots,n_{i})\otimes c^{{\mathcal F}_{i+1},\ldots,{\mathcal F}_{i+j}}(n_{i+1},\ldots,n_{i+j})),\nonumber
\end{align}
\begin{align}
&(c^{i}\circ (c^{j_{1}},\ldots,c^{j_{r}}))^{{\mathcal F}_{1},\ldots,{\mathcal F}_{i+j_{1}+\cdots +j_{r}-r}}(n_{1},\ldots,n_{i+j_{1}+\cdots +j_{r}-r})\label{brac}\\
&:= \sum_{0\leq k_{1}\leq \ldots\leq k_{r}\leq i+j_{1}+\cdots +j_{r}-r}(-1)^{\sum_{p=1}^{r}(j_{p}-1)k_{p}}\nonumber\\
& c^{{\mathcal F}_{1},\ldots,{\mathcal F}_{k_{1}},{\mathcal F}_{k_{1}+1}\otimes\cdots\otimes{\mathcal F}_{k_{1}+j_{1}},\ldots,{\mathcal F}_{k_{r}},{\mathcal F}_{k_{r}+1}\otimes\cdots\otimes{\mathcal F}_{k_{r}+j_{r}},\ldots,{\mathcal F}_{i+j_{1}+\cdots +j_{r}-l}}\nonumber\\
&(n_{1},\ldots,n_{k_{1}}, c^{{\mathcal F}_{k_{1}+1},\ldots,{\mathcal F}_{k_{1}+j_{1}}}(n_{k_{1}+1},\ldots,n_{k_{1}+j_{1}}),\nonumber\\
& \ \ \ \ \ \ldots,n_{k_{r}},c^{{\mathcal F}_{k_{r}+1},\ldots,{\mathcal F}_{k_{r}+j_{r}}}(n_{k_{r}+1},\ldots,n_{k_{r}+j_{r}}),\ldots, n_{i+j_{1}+\cdots +j_{r}-r}).\nonumber
\end{align}

As in the classical case of the Hochschild complex of an associative algebra \cite{ger}\cite{ger-vor} we see that the following crucial relations hold.
\begin{lemma}The  monoidal functor structure (\ref{mon1})-(\ref{mon2}), the Hochschild type differential (\ref{diff1})-(\ref{diff2}) and the Gerstenhaber-Voronov type operations (\ref{prod})-(\ref{brac}) on the cochains (\ref{cochain1})-(\ref{cochain2}) satisfy the following identities
\begin{align}
c^{i}c^{j} & = (-1)^{i}\mu\circ (c^{i}, c^{j}),\\
\mathrm{d}c^{k} & = \mu\circ c^{k}-(-1)^{k-1}c^{k}\circ \mu.
\end{align}
Moreover, the Gerstenhaber  type operation (\ref{brac}) (the substitution $\circ $) satisfies the pre-Jacobi identity
\begin{align}
&(c^{i}\circ (c^{j_{1}},\ldots,c^{j_{r}}))\circ (c^{k_{1}},\ldots,c^{k_{s}}) \\
&:= \sum_{0\leq s_{1}\leq \ldots\leq s_{r}\leq s}
(-1)^{\sum_{p=1}^{r}(j_{p}-1)\sum_{q=1}^{q_{p}}(k_{q}-1)}\nonumber\\
& c^{i}\circ (c^{k_{1}},\ldots,c^{k_{q_{1}}}, c^{j_{1}}\circ (c^{k_{q_{1}+1}},\ldots),\nonumber\\
& \ \ \ \ \ \ \ \ \ \ \ \ldots,c^{k_{q_{r}}}, c^{j_{r}}\circ (c^{k_{q_{r}+1}},\ldots),\ldots, c^{k_{s}}).\nonumber
\end{align}
\end{lemma}
Therefore the same arguments as in \cite{ger-vor} prove the following theorem.
\begin{theorem}
The
Hochschild complex of a monoidal functor with the product (\ref{prod}) and substitutions (\ref{brac}) is a strong homotopy Gerstenhaber algebra.
\end{theorem}
In particular, as for the strong homotopy Gerstenhaber algebra of an associative algebra \cite{ger-vor}, this implies the following corollary \cite{Yett,Yett'}
\begin{corollary} The product (\ref{prod}) and the bracket $\{ c^{i}, c^{j}\} := c^{i}\circ c^{j} - (-1)^{(i+1)(j+1)}c^{j}\circ c^{i}$ define the structure of a Gerstenhaber
algebra on the Hochschild cohomology of a monoidal functor.
\end{corollary}

\subsection{Hochschild cohomology and deformations}
Exactly as in the  case of associative algebras \cite{Ger} the following results relating  infinitesimal deformations and Hochschild cohomology still hold \cite{Yett,Yett'}.
\begin{theorem}
There is a one-to-one correspondence between the group of equivalence
classes of 1-st infinitesimal deformations of a given monoidal functor and its second Hochschild cohomology group.
\end{theorem}

Moreover, succesive lifting of equivalences ``modulo $t^{2}$, $t^{3}$, ...'' of formal deformations \cite{Ger} obviously makes perfect sense also for monoidal functors and the following theorem still holds.

\begin{theorem}
If the second Hochschild cohomology of a given monoidal functor vanishes then its all formal
deformations are equivalent.
\end{theorem}

Assume now that we have an $n$-th infinitesimal deformation $(\mu_{0}, \ldots, \mu_{n})$ which lifts to an $(n+1)$-st infinitesimal deformation $(\mu_{0}, \ldots, \mu_{n+1})$. Then the identity (\ref{defass}) for $(\mu_{0}, \ldots, \mu_{n+1})$ can be rewritten as
\begin{align}\sum_{j=1}^{n}\mu_{n+1-j}^{{\mathcal F}_{1},{\mathcal F}_{2}\otimes{\mathcal F}_{3}}
(Id^{f_{*}{\mathcal F}_{1}}\otimes \mu_{j}^{{\mathcal F}_{2},{\mathcal F}_{3}})
-\mu_{n+1-j}^{{\mathcal F}_{1}\otimes{\mathcal F}_{2},{\mathcal F}_{3}}
(\mu_{j}^{{\mathcal F}_{1},{\mathcal F}_{2}}
\otimes Id^{f_{*}{\mathcal F}_{3}})=(d\mu_{n+1})^{{\mathcal F}_{1},{\mathcal F}_{2},{\mathcal F}_{3}}.\label{extdefass}
\end{align}
As in the classical case \cite{Ger}, by a straightforward computation, one can check that the left hand side of (\ref{extdefass}) is a 3-cocycle in the Hochschild complex of a given monoidal functor and hence defines a 3-dimensional cohomology class $\mathfrak{O}_{n}$. Therefore (\ref{extdefass}) reads as follows.
\begin{theorem}
$\mathfrak{O}_{n}$ is an obstruction for lifting a given n-th infinitesimal deformation to some $(n+1)$-st one.
\end{theorem}
\subsection{Comparison with the classical Hochschild complex} Since a monoidal functor transforms monoids into monoids the evaluation of $f_{*}$ on the monoidal unit ${\mathcal O}_{X}$ is a monoid with the unit and multiplication
\begin{align}
\eta: &\ \  {\mathcal O}_{Y}\rightarrow f_{*}{\mathcal O}_{X},\\
\mu^{{\mathcal O}_{X}, {\mathcal O}_{X}}: & \ \ f_{*}{\mathcal O}_{X}\otimes f_{*}{\mathcal O}_{X}\rightarrow f_{*}({\mathcal O}_{X}\otimes {\mathcal O}_{X})=f_{*}{\mathcal O}_{X}.
\end{align} It has its own Hochschild cochain complex with its canonical strong homotopy Gerstenhaber algebra structure. It is easy to check that the  two strong homotopy Gerstenhaber algebra structures are compatible.
\begin{proposition}
The map from the Hochschild cochain complex of a given monoidal functor $f_{*}$ to Hochschild cochain complex of a monoid $f_{*}{\mathcal O}_{X}$ given by the evaluation of natural transformations at the monoidal unit
\begin{align}
c^{k}\mapsto c^{{\mathcal O}_{X},\ldots, {\mathcal O}_{X}}
\end{align}
is a homomorphism of strong homotopy Gerstenhaber algebras.   
\end{proposition}
In particular, the usual Hochschild complex of the monoid $f_{*}{\mathcal O}_{X}$ is only a one particular component of the Hochschild complex of the monoidal functor $f_{*}$.
\subsection{Poisson functors} Let categories ${\rm Qcoh}_{X}$, ${\rm Qcoh}_{Y}$ be symmetric with symmetries
\begin{align}\sigma_{X}^{{\mathcal F}_{1},{\mathcal F}_{2}}: {\mathcal F}_{1}\otimes{\mathcal F}_{2}\rightarrow {\mathcal F}_{2}\otimes{\mathcal F}_{1},\ \ \sigma_{Y}^{{\mathcal G}_{1},{\mathcal G}_{2}}: {\mathcal G}_{1}\otimes{\mathcal G}_{2}\rightarrow {\mathcal G}_{2}\otimes{\mathcal G}_{1} \label{symm}\end{align}

The symmetry of the monoidal functor of $f_{*}: {\rm Qcoh}_{X}\rightarrow
{\rm Qcoh}_{Y}$ is described by the identity
\begin{align}f_{*}(\sigma_{X}^{{\mathcal F}_{1},{\mathcal F}_{2}})\mu^{{\mathcal F}_{1},{\mathcal F}_{2}}=\mu^{{\mathcal F}_{2},{\mathcal F}_{1}}\sigma_{Y}^{f_{*}{\mathcal F}_{1},f_{*}{\mathcal F}_{2}}.\label{symmon}\end{align}
Note that if the monoidal functor $f_{*}$ is symmetric then the  monoid $f_{*}{\mathcal O}_{X}$ is commutative.

\vspace{3mm}
\paragraph{\textbf{Definition}} We say that a symmetric monoidal functor $f_{*}$ is \emph{Poisson} if there is given a natural transformation $\pi$ of bifunctors
\begin{align}\pi^{{\mathcal F}_{1},{\mathcal F}_{2}}: f_{*}{\mathcal
F}_{1}\otimes f_{*}{\mathcal F}_{2}\rightarrow f_{*}({\mathcal
F}_{1}\otimes {\mathcal F}_{2})\label{pois}\end{align}
satisfying the following identities:

\emph{(skew symmetry)}
\begin{align}\pi^{{\mathcal F}_{1},{\mathcal F}_{2}}=-f_{*}(\sigma_{X}^{{\mathcal F}_{2},{\mathcal F}_{1}})\pi^{{\mathcal F}_{2},{\mathcal F}_{1}}\sigma_{Y}^{f_{*}{\mathcal F}_{1},f_{*}{\mathcal F}_{2}},\label{skew}\end{align}

\emph{(Jacobi identity)}
\begin{align}\pi^{{\mathcal F}_{1},{\mathcal F}_{2}\otimes {\mathcal F}_{3}}(Id^{f_{*}{\mathcal F}_{1}}\otimes \pi^{{\mathcal F}_{2},{\mathcal F}_{3}})-\pi^{{\mathcal F}_{1}\otimes{\mathcal F}_{2}, {\mathcal F}_{3}}(\pi^{{\mathcal F}_{1},{\mathcal F}_{2}}\otimes Id^{f_{*}{\mathcal F}_{3}})\label{jac}\\
=f_{*}(\sigma_{X}^{{\mathcal F}_{2},{\mathcal F}_{1}}\otimes Id^{{\mathcal F}_{3}})\pi^{{\mathcal F}_{2},{\mathcal F}_{1}\otimes {\mathcal F}_{3}}(Id^{f_{*}{\mathcal F}_{2}}\otimes \pi^{{\mathcal F}_{1}, {\mathcal F}_{3}})(\sigma_{Y}^{f_{*}{\mathcal F}_{1},f_{*}{\mathcal F}_{2}}\otimes Id^{f_{*}{\mathcal F}_{3}}),\nonumber\end{align}

\emph{(derivation)}
\begin{align}\pi^{{\mathcal F}_{1},{\mathcal F}_{2}\otimes {\mathcal F}_{3}}(Id^{f_{*}{\mathcal F}_{1}}\otimes \mu^{{\mathcal F}_{2},{\mathcal F}_{3}})-\mu^{{\mathcal F}_{1}\otimes{\mathcal F}_{2}, {\mathcal F}_{3}}(\pi^{{\mathcal F}_{1},{\mathcal F}_{2}}\otimes Id^{f_{*}{\mathcal F}_{3}})\label{der}\\
=f_{*}(\sigma_{X}^{{\mathcal F}_{2},{\mathcal F}_{1}}\otimes Id^{{\mathcal F}_{3}})\mu^{{\mathcal F}_{2},{\mathcal F}_{1}\otimes {\mathcal F}_{3}}(Id^{f_{*}{\mathcal F}_{2}}\otimes \pi^{{\mathcal F}_{1}, {\mathcal F}_{3}})(\sigma_{Y}^{f_{*}{\mathcal F}_{1},f_{*}{\mathcal F}_{2}}\otimes Id^{f_{*}{\mathcal F}_{3}}).\nonumber\end{align}

\begin{proposition}
If the symmetric monoidal functor $f_{*}$ is Poisson then the commutative monoid $f_{*}{\mathcal O}_{X}$ is Poisson as well, with the Poisson bracket
\begin{align}
\pi^{{\mathcal O}_{X}, {\mathcal O}_{X}}: & \ \ f_{*}{\mathcal O}_{X}\otimes f_{*}{\mathcal O}_{X}\rightarrow f_{*}({\mathcal O}_{X}\otimes {\mathcal O}_{X})=f_{*}{\mathcal O}_{X}.
\end{align}
\end{proposition}

\subsection{Quasi-classical limit} A similar calculation as in the case of a formal associative deformation of a commutative algebra one can easily check the following fact.

\begin{theorem}
Given a formal deformation $(\mu_{0}, \mu_{1}, \ldots )$ of a symmetric monoidal structure $(\eta, \mu)$ on $f_{*}$ the natural transformation $\pi$ of bifunctors
\begin{align}
\pi^{{\mathcal F}_{1}, {\mathcal F}_{2}}:= \mu_{1}^{{\mathcal F}_{1}, {\mathcal F}_{2}}-f_{*}(\sigma_{X}^{{\mathcal F}_{2},{\mathcal F}_{1}})\mu_{1}^{{\mathcal F}_{2}, {\mathcal F}_{1}}\sigma_{Y}^{f_{*}{\mathcal F}_{1},f_{*}{\mathcal F}_{2}}
\end{align}
makes $f_{*}$ a Poisson functor. 
\end{theorem}

\end{document}